\documentclass{elsart3-1}

\usepackage{amssymb,amsmath}

\usepackage[english,francais]{babel}

\newtheorem{theorem}{Theorem}[section]
\newtheorem{lemma}[theorem]{Lemma}
\newtheorem{e-proposition}[theorem]{Proposition}

\newtheorem{e-definition}[theorem]{Definition\rm}
\newtheorem{remark}{\it Remark\/}

\renewcommand{\ge}{\varepsilon}
\newcommand{\R}{\mathbb{R}}
\newcommand{\Rn}{{\mathbb{R}^N}}
\newcommand{\D}{{\mathcal{D}}}
\newcommand{\subD}{{\mathcal{D}}}

\def\og{\leavevmode\raise.3ex\hbox{$\scriptscriptstyle\langle\!\langle$~}}
\def\fg{\leavevmode\raise.3ex\hbox{~$\!\scriptscriptstyle\,\rangle\!\rangle$}}

\begin{document}

\begin{frontmatter}


\selectlanguage{english}
\title{Remarks on a Hardy-Sobolev inequality}




\selectlanguage{english}
\author[Simone Secchi]{Simone Secchi}
\ead{secchi@mail.dm.unipi.it}
\author[Didier Smets]{Didier Smets}
\ead{smets@ann.jussieu.fr}
\author[Michel Willem]{Michel Willem}
\ead{willem@math.ucl.ac.be}

\address[Simone Secchi]{Dipartimento di Matematica, Universit\`a di
Pisa, via F. Buonarroti, 2 I-56127 Pisa, Italy}
\address[Didier Smets]{Laboratoire J.-L. Lions, Universit\'e Paris 6,
4 pl. Jussieu BC 187, 75252 Paris Cedex 05, France} 
\address[Michel Willem]{Universit\'e Catholique de Louvain, 2 chemin
du cyclotron, 1348 Louvain-la-Neuve, Belgium} 

\begin{abstract}
We compute the optimal constant for a generalized Hardy-Sobolev 
inequality, and using the product of two symmetrizations
we present an elementary proof of the symmetries of some optimal
functions. This inequality was motivated by a nonlinear elliptic
equation arising in astrophysics.

\vskip 0.5\baselineskip

\noindent{\large Remarques sur une in\'egalit\'e de Hardy-Sobolev}
\vskip 0.5\baselineskip

\selectlanguage{francais}
\noindent{\bf R\'esum\'e}
\vskip 0.5\baselineskip
\noindent
Nous calculons la meilleure constante dans une in\'egalit\'e de
Hardy-Sobolev g\'en\'eralis\'{e}e, et en utilisant le produit de deux
sym\'{e}trisations, nous montrons de mani\`ere \'el\'ementaire la
sym\'{e}trie de certaines fonctions optimales. Cette in\'{e}galit\'{e}
est motiv\'{e}e par une \'{e}quation elliptique non-lin\'{e}aire en astrophysique. 
\end{abstract}
\end{frontmatter}

\selectlanguage{francais}
\section*{Version fran\c{c}aise abr\'eg\'ee}
Nous d\'eterminons la constante optimale $C$
dans l'in\'egalit\'e de type Hardy-Sobolev 
\begin{equation}
\label{eq:1fr}
\int_\Rn \frac{|u(x)|^q}{|y|^\beta}\, dx \leq C \int_\Rn |\nabla
u|^p\,  dx,
\end{equation}
o\`u $x=(y,z)\in\R^k\times\R^{N-k}$. Celle-ci est donn\'ee par la valeur
\[
C=\frac{p^p}{(k+p)^p}.
\]
Nous consid\'erons \'egalement la sym\'etrie de certaines fonctions
optimales. En utilisant le produit de deux sym\'etrisations, nous
prouvons qu'elles ne d\'ependent que de $(|y|, |z|)$.

\selectlanguage{english}
\section{Introduction}
\label{}
\label{intro}
The Hardy-Sobolev inequality
\begin{equation}
  \label{eq:SH}
  \int_\Rn \frac{|u(x)|^q}{|y|^\beta}\, dx \leq C \int_\Rn |\nabla
u|^p\,  dx
\end{equation}
where $x=(y,z)\in\R^k\times\R^{N-k}$ was studied by Badiale and
 Tarantello in \cite{BT}.

Our aim is to solve two open problems contained in \cite{BT}. First we
 compute the optimal value of the constant $C$ in Equation \eqref{eq:SH}
 in the case of Hardy's inequality, namely $p=q=\beta$.
In fact we prove a more general inequality with optimal constant in
 Section 2.
Second, in Section 3, we consider the symmetry of the optimal functions. Using
the ``product'' of two symmetrizations, we prove that some optimal
functions depend only on $(|y|,|z|)$. 

\begin{center}
{\bf Notation}
\end{center}
\begin{description}
\item{\textbullet} When we write an integral like $\int_\Rn |u|$, we
 mean that the integral is taken with respect to the Lebesgue measure on
 $\Rn$. 
\item{\textbullet} $\D(\Rn)$ denotes the space of test functions, namely
\[
\D(\Rn)=\{ u\in C^\infty (\Rn) \mid \operatorname{supp} u \text{ is
 compact}\}.
\]
\item{\textbullet} $D^{1,p}(\Rn)$ is the closure of $\D(\Rn)$ with
 respect to the norm
\[
\|u\|_{D^{1,p}}=\left( \int_\Rn |\nabla u|^p  \right)^{1/p}.
\]
\item{\textbullet} $W^{1,p}(\Rn)$ is the usual Sobolev space, namely
 the closure of $\D(\Rn)$ with respect to the norm
\[
\|u\|_{W^{1,p}}=\left( \int_\Rn |\nabla u|^p  +\int_\Rn |u|^p
 \right)^{1/p}.
\]
\end{description}

\section{Generalized Hardy inequality}

If $1\leq k\leq N$, we will write a generic point $x\in\Rn$ as
 $x=(y,z)\in \R^k\times\R^{N-k}$.

\medskip

\begin{theorem}\label{th:1}
Let $1<p<\infty$ and $\alpha +k>0$. Then, for each $u\in \D(\Rn)$ the
following inequality holds: 
\begin{equation}
  \label{eq:hardy}
  \int_\Rn |u(x)|^p |y|^\alpha \, dx \leq \frac{p^p}{(\alpha+k)^p}
 \int_\Rn |\nabla u(x)|^p |y|^{\alpha+p}\, dx.
\end{equation}
Moreover, the constant $\frac{p^p}{(\alpha+k)^p}$ is optimal.
\end{theorem}

\medskip

\begin{lemma}\label{lem:1}
Inequality \eqref{eq:hardy} holds when $k=N$.
\end{lemma}

\noindent{\it Proof : }
See \cite{K,W} and also \cite{M}, where many generalizations are proved.
\hfill\mbox{}\qed

\medskip

\begin{lemma}\label{lem:2}
When $k=N$ the constant $\frac{p^p}{(\alpha+N)^p}$ in inequality
 \eqref{eq:hardy} is optimal.
\end{lemma}

\noindent{\it Proof : }
Consider the family of functions
\begin{equation*}
u_\ge (x)=
\begin{cases}
1 &\text{ if }|x|\leq 1,\\
|x|^{-\frac{\alpha+N}{p}-\ge} &\text{ if }|x|>1.
\end{cases}
\end{equation*}
and pass to the limit as $\ge\to 0$.
\hfill\mbox{}\qed

\medskip

\begin{lemma}\label{lem:3}
Inequality \eqref{eq:hardy} holds for any $1\leq k\leq N$.
\end{lemma}

\noindent{\it Proof : }
For every $u\in \D(\Rn)$ we have by Lemma \ref{lem:1}
\begin{eqnarray*}
  \int_\Rn |u(x)|^p |y|^\alpha\, dx &=& \int_{\R^{N-k}}\, dz \int_{\R^k}
 |u(x)|^p |y|^\alpha\, dy\\
&\leq& \frac{p^p}{(\alpha+k)^p} \int_{\R^{N-k}}\, dz \int_{\R^k} |\nabla_y
 u(x)|^p |y|^{\alpha+p}\, dy \\
&\leq& \frac{p^p}{(\alpha+k)^p} \int_\Rn |\nabla u(x)|^p |y|^{\alpha+p}\,
 dy.
\end{eqnarray*}
\hfill\mbox{}\qed

\medskip

\begin{lemma}\label{lem:4}
The constant $\frac{p^p}{(\alpha+k)^p}$ in inequality \eqref{eq:hardy} 
is optimal.
\end{lemma}

\noindent{\it Proof : }
Let us choose $u\colon (y,z)\mapsto v(y)w(z)$ with $v\in \D(\R^k)$ and
 $w\in \D(\R^{N-k})$. It is clear that
\[
\int_\Rn |\nabla u(x)|^p |y|^{\alpha+p}\, dx = \int_\Rn \left( |\nabla
 v(y)|^2 w(z)^2 + |\nabla w(z)|^2 v(y)^2  \right)^\frac{p}{2}
 |y|^{\alpha+p}\, dx,
\]
and
\[
\int_\Rn |u(x)|^p|y|^\alpha\, dx = \int_{\R^k} |v(y)|^p |y|^\alpha \, dy
 \int_{\R^{N-k}}|w(z)|^p \, dz.
\]
If we consider the convex function
\begin{align*}
F\colon [0,+\infty) \times [0,+\infty) &\longrightarrow [0,+\infty)\\
(s,t) &\mapsto \left( s^2+t^2    \right)^\frac{p}{2},
\end{align*}
we get that
\begin{equation}
  \label{eq:conv}
  \left( s^2+t^2    \right)^\frac{p}{2}\leq (1-\lambda)^{1-p}s^p +
 \lambda^{1-p}t^p,
\end{equation}
for all $s,t\geq 0$ and $0<\lambda <1$.
Hence we obtain, for $0<\lambda <1$,
\begin{equation*}
\frac{\int_\Rn |\nabla u|^p |y|^{\alpha+p}\, dx}{\int_\Rn |u|^p
 |y|^\alpha \, dx}\leq (1-\lambda)^p \frac{\int_{\R^k}|\nabla v|^p
 |y|^{\alpha+p} \, dy}{\int_{\R^k}|v|^p |y|^\alpha}
+ \lambda^p \frac{\int_{\R^{N-k}}|\nabla w|^p\,
 dz}{\int_{\R^{N-k}}|w|^p\, dz}\frac{\int_{\R^k} |v|^p|y|^{\alpha+p}\,
 dy}{\int_{\R^k}|v|^p |y|^\alpha\, dy}.
\end{equation*}
Since
\[
\inf_{\substack{w\in \subD(\R^{N-k})\\w\neq 0}}
 \frac{\int_{\R^{N-k}}|\nabla w|^p\, dz}{\int_{\R^{N-k}}|w|^p\, dz} = 0,
\]
we obtain, for $0<\lambda <1$,
\[
\inf_{\substack{u\in \subD(\R^{N})\\u\neq 0}} \frac{\int_\Rn |\nabla 
u|^p
 |y|^{\alpha+p}\, dx}{\int_\Rn |u|^p |y|^\alpha \, dx} \leq
 (1-\lambda)^p \inf_{\substack{v\in \subD(\R^{k})\\v\neq 0}}
 \frac{\int_{\R^k}|\nabla v|^p |y|^{\alpha+p} \, dy}{\int_{\R^k}|v|^p
 |y|^\alpha}.
\]
By letting $\lambda \to 0$, we deduce from Lemma 2 that
\[
\inf_{\substack{u\in \D(\R^{N})\\u\neq 0}} \frac{\int_\Rn |\nabla u|^p
 |y|^{\alpha+p}\, dx}{\int_\Rn |u|^p |y|^\alpha \, dx} \leq
 \frac{(\alpha+k)^p}{p^p}.
\]
The conclusion then follows from Lemma \ref{lem:3}.
\hfill\mbox{}\qed

\medskip 

\begin{remark}{\rm
When $k=N$ and $\alpha=-p$, inequality \eqref{eq:SH} is the classical
 Hardy inequality (see \cite{K} or \cite{W}).
}\end{remark}

\begin{remark}{\rm
When $2\leq k\leq N$ and $\alpha=-p$, inequality \eqref{eq:SH} was
 conjectured by Badiale and Tarantello in \cite{BT}.
}\end{remark}

\begin{remark}{\rm
This method is applicable in many other problems. A simple example is
 that for any open subset $\Omega\subset \R^M$,
\[
\inf_{\substack{u\in \subD(\Omega\times \Rn)\\
 \int\limits_{\Omega\times\Rn}|u|^p=1}} \int_{\Omega\times\Rn} |\nabla 
u|^p =
 \inf_{\substack{v\in \subD(\Omega)\\ \int\limits_{\Omega}|v|^p=1}} 
\int_\Omega
 |\nabla v|^p.
\]
}\end{remark}

\section{Cylindrical symmetry}

In this section, we consider the minimization problem
\begin{equation}
  \label{eq:S}
  S=S(N,p,k,\beta)=\inf \left\{ \int_\Rn |\nabla u|^p \mid u\in
 D^{1,p}(\Rn) \quad \text{and}\,\,\, \int_\Rn \frac{|u|^q}{|y|^\beta}=1
\right\}
\end{equation}
where

\vspace{10pt}

\noindent \textbf{(H)} \qquad\qquad\quad $0\leq \beta <k,\quad\beta \leq
 p, \quad q=q(N,p,\beta)=\frac{p(N-\beta)}{N-p}$.

\vspace{10pt}

Let $u\in L^1(\Rn)$ be a non-negative function. Let us denote by
 $u^\star(\cdot,z)$ the Schwarz symmetrization of $u(\cdot,z)$ and by
 $u^{\star\star}(y,\cdot)$ the Schwarz symmetrization of
 $u^\star(y,\cdot)$. It is clear that $u^{\star\star}$ depends only on
 $(|y|,|z|)$. Let us define
\[
D_{\star\star}^{1,p}(\Rn)=\left\{ u\in D^{1,p}(\Rn) \mid
 u=u^{\star\star}
\right\}
\]
and
\begin{equation}
S^{\star\star}=S^{\star\star}(N,p,k,\beta)=\inf \bigg\{ \int_\Rn
 |\nabla u|^p \mid u\in D_{\star\star}^{1,p}(\Rn)
\text{and}\,\,\, \int_\Rn \frac{|u|^q}{|y|^\beta}=1
\bigg\}.
\end{equation}

\medskip

\begin{theorem}\label{th:6}
Under assumption \textbf{(H)}, $S=S^{\star\star}$.
\end{theorem}

\medskip 

Since it is clear that $S\leq S^{\star\star}$, Theorem \ref{th:6} 
follows from a density argument and the next lemma.

\medskip

\begin{lemma}
Under assumption \textbf{(H)}, for each $u\in \D(\Rn)$, $u$ 
non-negative, then
\[
\int_\Rn |\nabla u^{\star\star}(x)|^p \leq \int_\Rn |\nabla u(x)|^p,
\]
\[
\int_\Rn \frac{|u^{\star\star}(x)|^q}{|y|^\beta}\geq \int_\Rn
 \frac{|u(x)|^q}{|y|^\beta}.
\]
\end{lemma}

\noindent{\it Proof : }
a) By the Polya-Szeg\"o inequality for Steiner symmetrization
 (see \cite{BS}, Theorem 8.2), $u^{\star}$ and $u^{\star\star}$ belong to
 $W^{1,p}(\Rn)$ and
\[
\int_\Rn |\nabla u^{\star\star}|^p \leq \int_\Rn |\nabla u^{\star}|^p
 \leq \int_\Rn |\nabla u|^p.
\]
b) Let $R>0$ be such that
\[
\operatorname{supp} u \subset B(0,R)\times \R^{N-k},
\]
and define
\[
v(y,z)=
\begin{cases}
\frac{1}{|y|^\beta} &\text{if } |z|\leq R,\\
0 &\text{otherwise}.
\end{cases}
\]
Since $v^{\star\star}=v$, it follows from the general Hardy-Littlewood
 inequality (see \cite{W}) that
\begin{equation*}
\int_\Rn \frac{|u^{\star\star}|^q}{|y|^\beta}=\int_\Rn \left( |u|^q 
 \right)^{\star\star} v^{\star\star}
\geq \int_\Rn |u|^q v = \int_\Rn \frac{|u|^q}{|y|^\beta}.
\end{equation*}
The proof is complete.\hfill\mbox{}\qed

\medskip

\begin{remark}{\rm
If $\beta <p$, then $S$ is achieved (see \cite{BT}).
}\end{remark}
\begin{remark}{\rm
It is proved in \cite{BT}, Theorem 5.3, that when $\beta<p=2$, the
 optimal function for $S$ satisfies $u=u^\star$.
}\end{remark}

\begin{remark}{\rm
The use of the "product"
of two symmetrizations is applicable to many other problems.
A simple example is given by constraints of the form
\[
\int_\Rn |u(y,z)|^q g(|y|)h(|z|)=1
\]
where $g$ and $h$ are non-increasing.
}\end{remark}

\noindent{\it Added in proof : } 
After completing this note, we learned of the preprint \cite{MS},
where the authors show that all minimizers for the
Hardy-Sobolev inequality must have the aforementioned symmetry, after
a possible translation in the $z$ variable.




\end{document}